\mag=\magstephalf
%\mag=\magstep1
\pageno=1
\input amstex
%\baselineskip = 1.0 true cm
\documentstyle{amsppt}
\TagsOnRight
\interlinepenalty=1000
%\hsize=6.5truein
%\voffset=0.3truein
%\hoffset=0.2truein
%\vsize =9.8truein
\NoRunningHeads
%\pagewidth{14cm}
%\pageheight{21.5cm}
%\vcorrection{-1.0cm}
%\hcorrection{-1.2cm}

\pagewidth{16.5 truecm}
\pageheight{23.0 truecm}
\vcorrection{-1.0cm}
\hcorrection{-0.5cm}
%\advance\vsize by -\voffset
%\advance\hsize by -\voffset
%\baselineskip = 1.0 true cm
\nologo

%Draft 99/09/30
\NoBlackBoxes

\def\tvskip{\vskip 0.5 cm}

\define \bee{{\bold e}}
\define \bEE{{\bold E}}
\define \tr{{\roman {tr}}}
%\define \det{{\roman det}}

\define \ind{{\roman {ind}}}

\define \dd{{\roman d}}

\define \EE{\Bbb E}
\define \RR{\Bbb R}
\define \ZZ{\Bbb Z}

\define \NN{\Bbb N}

\define \XiT{\Xi_{T_S}}

\define \CC{\Bbb C}

\define \conf{\roman{conf}}

\define \SS{\roman{S}}
\define \fg{\frak g}

\define \Hom{\roman{Hom}}

\define \SHom{\Cal{H}om}
\define \END{{\roman{END}}}

\define \CLIF{{\roman{CLIFF}}}
\define \CLIFs{{\text{\it{CLIFF}}}}
\define \Clifs{{\Cal{C}\text{\it{liff}}}}
\define \Spin{{\Cal{S}\text{\it{pin}}}}
\define \SpinC{{\Cal{S}\text{\it{pin}}^{\Bbb C}}}
\define \SPIN{{\roman{SPIN}}}
\define \SPINs{{\text{\it{SPIN}}}}
\define \SPINCs{{\text{\it{SPIN}}}^{\Bbb C}}
\define \EMV{{\roman{UNV}}}
\define \SPINC{{\roman{SPIN}^{\Bbb C}}}

\define \HH{{\roman{H}}}
\define \SO{\roman{SO}}
\define \OO{\roman{O}}
\define \SSSO{\Cal S \Cal O}
\define \SSO{\Cal O} 
\define \Spect{\roman{Spect}}

\define \SPAN{\roman{span}}
\define \LCconeM{{\nabla_M^{\roman{LC}}}}
\define \CLconeM{{\nabla_M^{\roman{Clf}}}}
\define \SPconeM{{\nabla_M^{\roman{Spin}}}}
\define \LCcone{\nabla^{\roman{LC}}}

\define \SPcone{\nabla^{\roman{Spin}}}
\define \SPHcone{\nabla^{\roman{SpinH}}}
\define \SPWcone{\nabla^{\roman{SpinW}}}
\define \Proof{ {\it Proof} }
\define \fvskip{\vskip 1.0 cm}
\define \CM{{\Cal C^\omega_M}}
\define \DM{{\Cal D_M}}
\define \DT{{\Cal D_{T_S}}}
\define \DS{{\Cal D_{S}}}

%\define \CMC{{{}^{\Bbb C}{\Cal C}^{\omega }_M}}
%\define \CUC{{{}^{\Bbb C}{\Cal C}^{\omega }_U}}
%\define \CSC{{{}^{\Bbb C}{\Cal C}^{\omega }_S}}
%\define \CTC{{{}^{\Bbb C}{\Cal C}^{\omega }_{T_S}}}

\define \CTS{{\Cal C^\omega_{T_S}}}
\define \CS{{\Cal C^\omega_S}}

\define \extp{{\roman{ext}}}
\define \intp{{\roman{int}}}
\define \cliffp{{\roman{cliff}}}
\define \cliffq{{\frak{c}}}

\define \Nabla{\ /\!\!\!\!\nabla}
\define \Partial{\ /\!\!\!\partial}

\define \DDM{{= \!\!\!\!\!\!\nabla}}
\define \Dir{\roman{Dir}}
%\define \DDM{{\diagdown \!\!\!\!\!\nabla}}
\define \dalpha{{\dot \alpha}}
\define \dbeta{{\dot \beta}}
\define \dgamma{{\dot \gamma}}

\define \tnablaM{{\nabla^{\roman{SA}}_M}}
\define \tnablaT{{\nabla^{\roman{SA}}_{T_S}}}
\define \tnablaS{{\nabla^{\roman{SA}}_{S}}}

\define \SGL{\Cal{G}\Cal{L}}
\define \DMod{{\Cal D}}
\define \III{{\Cal I}}
%\baselineskip = 1.0 true cm

\topmatter \title \endtitle {\centerline{\bf{ Submanifold Differential 
Operators  in $\DMod$-Module Theory II:}}

{\centerline{\bf{Generalized Weierstrass and 
Frenet-Serret Relations as Dirac Equations}}}

{}
{}
\vskip 0.5 cm
{\centerline{ S. MATSUTANI }

{\centerline{ 8-21-1 Higashi-Linkan Sagamihara 228-0811 JAPAN }

%\endaffil
\address 
\endaddress
%\centerline{\bf{Abstract}}
\abstract 

This article is one of a series of papers. For this decade, the Dirac 
operator on a submanifold has been studied as a restriction of the Dirac 
operator in $n$-dimensional euclidean space $\EE^n$ to a surface or a space
curve as physical models. These Dirac operators are identified with operators
of the Frenet-Serret relation for a space curve case and of the generalized
Weierstrass relation for a conformal surface case and completely represent
the submanifolds. For example, the analytic index of Dirac operator of a 
space curve is identified with its writhing number. As another example, the
operator determinants of the Dirac operators are closely related to 
invariances of the immersed objects, such as Euler-Bernoulli and Willmore 
functionals for a space curve and a conformal surface respectively. In this
article, we will give mathematical construction of the Dirac operator by 
means of $\Cal D$-module and reformulate my recent results mathematically.

\centerline{{\bf MCS Codes:} 32C38, 34L40, 35Q40}
%*32C38 Sheaves of differential operators and their modules
%*34L40 Particular operators (Dirac, one-dimensional Schr\"odinger, etc.) 
%*35Q40 Equations from quantum mechanics 
% 32C25 Analytic subsets and submanifolds 
% 32G10 Deformations of submanifolds and subspaces 
%*53B25 Local submanifolds [See also 53C40] 
% 53C12 Foliations (differential geometric aspects) [See also 57R30, 57R32] 
% 53C42 Immersions (minimal, prescribed curvature, tight, etc.) 
%[See also 49Q05, 49Q10, 53A10, 57R40, 57R42] 
% 58D10 Spaces of embeddings and immersions 

\centerline{{\bf Key Words:} Weierstrass Relation, Frenet-Serret Relation, 
Dirac operator, Submanifold, $\Cal D$-Module }

\endabstract

\endtopmatter
%\newpage
%\baselineskip = 0.8 true cm

\vskip 0.5 cm
{\centerline{\bf{\S 1. Introduction}}

\vskip 0.5 cm

In an earlier article in this series [I], we showed the construction of the
submanifold Schr\"odinger operator in terms of $\DMod$-module theory. In this
article, we will apply the scheme to the spin bundle to construct the Dirac
operator of submanifold. We will call the previous article [I] and its 
proposition or definition and  reference like (I-2.1) and [I-Mat2], which 
means proposition or definition 2-1 and reference [Mat2] in [I] respectively.

Applying the  quantum mechanical  scheme [I and its references] to Dirac 
operators for a restricted particle along a low-dimensional submanifold in
$n$-dimensional euclidean space $\EE^n$, we obtained natural Dirac operators
on curves in $\EE^n$ ($n\ge 2$) [Mat1-3, 5, 6, 8, 10, MT] and  on conformal
surfaces in $\EE^n$ ($n=3,4$) [BJ, Mat11, 14, 16]. In this decade, I have 
been studying these Dirac operators and investigating their properties. From
physical point of view, I showed that they exhibit the symmetry of 
corresponding submanifolds and found a non-trivial extension of Atiyah-Singer
type index theorem to submanifold [Mat2,6]. The Dirac operators of curves in
$\EE^n$ ($n\ge 2$) are related to the Frenet-Serret relations and are 
identified with the Lax operators of ($1+1$)-dimensional soliton equations,
{\it e.g.}, modified Korteweg-de Vries equation [MT, Mat2,8], nonlinear 
Schr\"odinger equation [Mat1, 3, 6], complex modified Korteweg-de Vries 
equation [Mat15], and so on [Mat3]. The Dirac operators on conformal surfaces
in $\EE^n$ ($n=3,4$) are concerned with the generalized Weierstrass equation
representing a surface [Mat11, 14,16] and also identified  with the Lax 
operators of  modified Novikov-Veselov (MNV) equations [Ko1, 2, T1, 2]. The
generalized Weierstrass equation is very interesting from the viewpoint of
immersion geometry and differential geometry [B, Ke, Ko1, 2, TK, T1, 2, PP]. 

For example, for a case that  $S$ is a conformal surface immersed in 
$\EE^3$, the generalized Weierstrass equation is given as the zero mode
of the Dirac operator,
$$
      \Nabla_{S \hookrightarrow \EE^3}^{\Dir} 
      =  \sqrt{-1}\rho^{1/2}\pmatrix  p & \partial \\
                 \overline \partial & p  \endpmatrix\rho^{-1/2},
$$                                           
where $(z,\bar z)$ is a complex parameterization of $S$, $p$ is 
$\rho^{1/2} H$ when the volume element of $S$  is given by dvol $= \rho
d z d \bar z$ and $H$ is the mean  curvature. This Dirac operator was
obtained by various methods,  {\it e.g.}, direct computations [Ko1, 2], 
quotation algebraic method [Pe, PP] spin bundle consideration [Fr]
and quantum mechanical scheme [BJ, Mat11,14,16].

Similarly we can obtain the generalized Weierstrass equation of $\EE^4$
case [Ko2, PP, Mat16]. The quantum mechanical scheme [Mat16] is one of the
easiest approaches. Thus in order that mathematicians  understand the
recent progress of quantum mechanical approaches, I believe that the 
construction scheme should be translated to mathematical language and
be more precise. 

In 1991, I applied the Jensen and Koppe [I-JK] and da Costa [I-dC]
 scheme of submanifold quantum mechanics to
the Dirac operator over a plane curve as a first example
and investigated it for this decade [MT, Mat2, 5, 10]. Here we will
summary those.
Along the scheme described in the introduction of [I], 
we obtain the Dirac operator on a real analytic curve in
$\EE^2$ as [MT, Mat2, 6, 9],
$$
   \Nabla_{C \hookrightarrow \EE^2}^{\Dir}=\pmatrix  v 
   &  -\sqrt{-1} \partial_{s} \\
                 -\sqrt{-1} \partial_{s} & v  \endpmatrix ,
$$
where $s$ is the arclength of the curve,
$v:=k/2$ and $k$ is the curvature of the curve. 
The Dirac equation as a zero mode of this Dirac operator,
$$
     \Nabla_{C \hookrightarrow \EE^2}^{\Dir} \psi = 0 ,
$$
is essentially  equivalent with the Frenet-Serret relation.
By considering the eigenvalue equations,
$$
     \Nabla_{C \hookrightarrow \EE^2}^{\Dir} \psi = \lambda \psi ,
$$
and deformation of the curve preserving $\lambda$, we obtain 
the MKdV equation,
$$
    \partial_t v + 6 v^2 \partial_1^3 v + \partial_1^3 v = 0,
$$
where $t$ is deformation parameter which can be physically
interpreted as the real time of the curve approximately [MT, Mat10] 
and an adiabatic parameter in the fermionic system [MT, Mat2, 8, 12, 17].

Furthermore operator determinant of 
$\Nabla_{C \hookrightarrow \EE^2}^{\Dir}$
is the invariant of the system, viz the Euler-Bernoulli functional
$$
   \log \det \Nabla_{C \hookrightarrow \EE^2}^{\Dir}
    = \frac{1}{2\pi} \int ds \ k^2(s),
$$
and its index is given as
$$
   \ind \Nabla_{C \hookrightarrow \EE^2}^{\Dir}
   = w_{C \hookrightarrow \EE^2},
$$
where $w_{C \hookrightarrow \EE^2}$ is the writhing number of the curve
[MT, Mat2, 6, 8, 14].
Moreover the moduli of the immersion curve is classified by
the Dirac operator [Mat12, 17]. 
Hence the interpretation of the Dirac operator means the investigation
of the curve itself. Hence I believe that my Dirac operator is very 
important in geometry. 

This procedure, as I conjectured in [Mat2], can be extended to a surface
case and we can obtain the generalized Weierstrass relation, 
Novikov-Veselov equation and Willmore functional. 

Here I will comment upon my philosophical idea on these studies. When I 
started these studies, I did not have any proper language to express my 
motivation. However now I can express my motivation why I started these.

In algebraic number theory, we encounter the Fermat problems 
whether there exists an integer pair $(x,y)\in \ZZ^2$ satisfied with
the relation
$x^2 + y^2 = p$ or $x^2 + 3 y^2 = p$ for a given prime number $p$ [Coh].
These are solved in quadratic form theory and studied by quadratic 
integer. For the case of $x^2 + y^2 = p$, if $p\equiv 1$ modulo $4 \ZZ$,
there exists such integer pair $(x,y)$. Due to the Legendre symbol 
$\left( \dfrac{ -1 }{p} \right) =1$ for $p\equiv 1$ modulo $4 \ZZ$, 
there is a number $z$ such that $z^2 +1 \equiv 0$ modulo $p \ZZ$. Hence
$p$ is not prime number if we consider it in $\ZZ[\sqrt{-1}]$; $p 
\ZZ[\sqrt{-1}]$ is not prime ideal in $\ZZ[\sqrt{-1}]$ and must be 
decomposed. (Since the equation $x^2 +y^2 \equiv 0$ modulo $p \ZZ$  is
decomposed to $(x+ \sqrt{-1} y)(x - \sqrt{-1} y) \equiv 0$ modulo 
$p\ZZ[\sqrt{-1}]$, it means that for $p\equiv 1$ modulo $4 \ZZ$, 
there exists $(x,y)\in \ZZ^2$ such that
$(x+ \sqrt{-1} y)(x - \sqrt{-1} y) \equiv 0$ modulo $p\ZZ[\sqrt{-1}]$.)
Similarly we will consider $x^2 + 3 y^2 = p$ as $(x + \sqrt{-3} 
y)(x-\sqrt{-3} y)\equiv 0$ modulo $p\ZZ[\sqrt{-3}]$. These Abel 
extensions decompose the set of prime ideals $\roman{Spec}(\ZZ)$, which
is a first step of class field theory.

Here we will note that $\sqrt{-1}$ and $\sqrt{-3}$ are generators of
orientation (charity or  complex conjugate) of 
$\ZZ[\sqrt{-1}]\subset \CC$ or  $\ZZ[\sqrt{-3}]\subset \CC$ and this 
extension $\ZZ[\sqrt{-1}]$ leads us to the  algebraic number theory, 
ideal theory and arithmetic geometrical theory. Further we will note 
that in algebraic number theory, sets of prime ideals (  {\it e.g.}, 
$\roman{Spec}(\ZZ)$, $\roman{Spec}(\ZZ[\sqrt{-1}])$, 
$\roman{Spec}(\ZZ[\sqrt{-3}])$, and so on) are more important than prime
numbers themselves.

Jensen and Koppe and da Costa [I-dC, I-JK] found the 
natural quadratic operator for a case of plane curve;
$$
   \Delta_{C \hookrightarrow \EE^2}:= \Delta_C +\frac{1}{2}k^2,
$$ 
where $\Delta_C$ is the Beltrami-Laplace operator of curve and
$\frac{1}{2}k^2$ is integrand of the Euler-Bernoulli functional.
The Euler-Bernoulli functional and $\Delta_C$ are invariant
of the system and hermite. 
Let the set of the spectrum of $\Delta_{C \hookrightarrow \EE^2}$
be denoted by $\Spect( \Delta_{C \hookrightarrow \EE^2})$. 
Then up to $\Spect( \Delta_{C \hookrightarrow \EE^2})$,
$\Delta_{C \hookrightarrow \EE^2}\equiv 0$ in a certain sense. 
By adding  a generator of orientation $k$ in $\EE^2$, 
($k(s)$ changes its sign
depending upon the orientation), we will encounter the  
Dirac operator $\Nabla_{C \hookrightarrow \EE^2}^{\Dir}$  [MT], 
which is not hermite,
$$
\frac{1}{2}\tr_{2 \times 2}((\Nabla_{C \hookrightarrow \EE^2}^{\Dir})^2)
=\Delta_{C \hookrightarrow \EE^2},
$$
 and 
$\Spect( \Delta_{C \hookrightarrow \EE^2})$ is 
quadratically decomposed to the spectrum
$\Spect( \Nabla_{C \hookrightarrow \EE^2}^{\Dir})$
 of $\Nabla_{C \hookrightarrow \EE^2}^{\Dir}$,
in a certain sense.
This is resemble with the quadratic form theory [Coh].

Of course, we might encounter "fail of unique factorization"
in this problem. As prime ideal must be more important than
prime number in algebraic number theory,
this problem should be also expressed by a module theory.
Thus I have chosen $\DMod$-module theory as mathematical expression
of my theory. Even though in noncommutative ring,
prime ideal becomes nonsense, I wish to construct a quadratic
form theory in this problem.

In other words, my motivation is arithmeticalization of quantum
physics or differential geometry and harmonic map theory.
I believe that the operator space of the  Dirac operator
plays the same role as $\ZZ[\sqrt{-1}]$ and $\ZZ[\sqrt{-3}]$.
This aspect can be regarded as another aspect of the program
of Eichler [E].

Recently number theory becomes gemetorized [E, Fa] and physics
and geometry also becomes arithmeticalized [Mat18 and references therein].
In future, both might be unified.
I hope that my approach might shed a new light on 
this way.

Contents are as follows.
\S 2 quickly reviews the Clifford algebra,
 the spin group and its spinor representation, which are
 well-established [ABS, BGV, G, Y]. The purpose of \S2 is
 to show notations in this article. 
In \S3, using well-known results on the Clifford, spin and
spinor bundle [ABS, BGV, G, Y], we will follow the argument
of Mallios [Mal1, 2] and formally construct
sheaves of spin, Clifford module and so on.
In \S 4, we will define the Dirac operator in a
submanifold and give theorems.

\centerline{\bf Acknowledgment}
\tvskip 
I will sincerely thank Prof.~B.~L.~Konopelchenko 
and Prof.~I.A.Taimanov for sending their works and
many encouragements 
and to Prof. F.~Pedit for giving a lecture on his quotation
theory and sending me the paper [PP].

I am grateful to Prof. K.~Tamano, Prof. Y.~\^Onishi,
Prof.~S.~Saito, Prof.T.~Tokihiro, Dr.H.~Tsuru, Prof.A.~Suzuki,
Prof.S.~Takagi, Prof.~H.~Kuratsuji, Prof.~K.~Sogo,
W.~Kawase and H.~Mitsuhashi for helpful discussions and comments in 
early stage of this study.

I thank Prof.~Y.~Ohnita, Prof.~M.~Guest, Dr.~R.~Aiyama and 
Prof.~K.~Akutagawa for inviting me their seminars and for critical 
discussions.

\vskip 0.5 cm
%\newpage

{\centerline{{\bf \S 2. Spinor Group}}
\tvskip

Let us review the spinor group [BGV, Ch, G].
 
\subheading{Definition 2.1}
{\it
\roster

Let $V$ be a real vector with a positive definite inner product $(,)$.

\item
We will define an equivalent relation $\sim_\CLIF$ by
$v \cdot v + (v,v) = 0 \ \text{for } v \in V$.

\item
The universal algebra generated by $V$ with $\RR$ coefficient
 is denoted by $\EMV(V)$.

\item
The Clifford algebra $\CLIF(V)$ is defined by
$$
\CLIF(V)=\EMV(V)/\sim_\CLIF .
$$
\endroster
}

\subheading{Lemma 2.2}

{\it
The relation $v \cdot v + (v,v) = 0 $
is given by the basis
$$
   e_i e_j + e_j e_i = -2\delta_{i j },
$$
where $\delta_{ i j}$ is Kronecker delta symbol.
}

\demo{Proof}{
First we will consider the case $v=e_i$. Then
we have $e_i \cdot e_i = -1$.
Next let $v= \sum v^i e_i$. The relation is
$v\cdot v=\sum_{i,j}v^i v^j e_i\cdot e_j +$
$\sum_i (v^i)^2=0$ and reduces to
$e_i e_j + e_j e_i=0$ for $i \neq j$.
}\enddemo
\tvskip

Let we will show the representation of $\CLIF(V)$
using endomorphism of the exterior algebra 
$\Lambda V$ of $V$, $\END(\Lambda V)$.

\subheading{Definition 2.3}
{\it
\roster
Let an orthonormal basis of $V$ be $\{e_j\}_{i=1,\cdots,n}$. 

\item
Let $\extp:V \to \END(\Lambda V)$ be exterior multiplication
on the left. Its action is 
given as,
$$
   \extp(e_1) (e_{i_1}\wedge \cdots
   \wedge e_{i_p})=\left\{
   \matrix 
   e_1\wedge e_{i_1}\wedge \cdots\wedge e_{i_p}, 
   & \text{if } i_1 \neq 1 \\
   0 , & \text{otherwise} 
   \endmatrix \right. .
$$

\item
Let $\intp:V \to \END(\Lambda V)$ 
be interior multiplication through the adjoint map
$*: V \to V^*:=\Hom(V, \RR)$.
Its action is 
given as
$$
   \intp(e_1) (e_{i_1}\wedge \cdots
   \wedge e_{i_p})=\left\{
   \matrix 
   e_{i_2}\wedge \cdots\wedge e_{i_p}, & \text{if } i_1 = 1 \\
   0 , & \text{otherwise} 
   \endmatrix\right. .
$$

\item
Let $\cliffp:V \to \END(\Lambda V)$ be clifford multiplication,
$$
    \cliffp(v) := \extp(v)-\intp(v).
$$
\endroster}

\tvskip

\subheading{Lemma 2.4}
{\it
\roster

\item
There is a ring endomorphism  by extension of
the action $\cliffp:V \to \END(\Lambda V)$
$$
   \cliffp : \CLIF(V) \to \END(\Lambda V).
$$

\item
$\Lambda V$ is $\CLIF(V)$-module.

\item There is a vector space isomorphism as $\CLIF(V)$-module, 
$$
   \frak{cliff} : \CLIF(V) \to\Lambda V, \quad
   (a \mapsto  \frak{cliff}(a):=\cliffp(a)(1)).
$$
Let us express $\cliffq :=(\frak{cliff})^{-1}$. 

\endroster
}

\demo{\Proof}
Noting $\extp(v)^2=\intp(v)^2\equiv 0$ for $v \in V$, we have a relation,
$$
   (\cliffp(v))^2 = -(\extp(v)\intp(v)+\intp(v)\extp(v))= -(v,v) I.
$$
First we will consider $\cliffp:\EMV(V) \to \END(V)$ as
for $w_1, w_2 \in \EMV(V)$, 
$\cliffp( w_1 + w_2)=\cliffp( w_1) +\cliffp(  w_2)$.
and $w_1, w_2 \in \EMV(V)$, $\cliffp( w v ) := \cliffp(w) \cliffp(v)$.
Then above relation is consist with $\sim_\CLIF$ 
and $\cliffp:\CLIF(V) \to \END(\Lambda V)$
is well-defined. Others are easy.
\qed \enddemo

\subheading{Definition 2.5}
{\it

We will decompose $\CLIF(V)$ as a module to submodule $\CLIF_+(V)$ 
consisting of multiplication of even elements of $V$ and submodule 
$\CLIF_-(V)$ consisting of multiplication of odd elements of $V$;
$$
   \CLIF(V)=\CLIF_+(V)\oplus\CLIF_-(V).
$$

}

\subheading{Definition 2.6 (Spin Group)}
{\it

Let $*$ be involution of $\CLIF(V) \to \CLIF(V)$
as $(e_{i_1} \cdots e_{i_p})^* =(e_{i_p} \cdots e_{i_1})$.

We will  define the spin group,
$$
   \SPIN(V) = \{ g \in \CLIF_+(V)\ |\ g^* g
   =1,\  g^* v g \in V,\  \forall v \in V\}.
$$
For $V=\RR^n$ case, we will denote it $\SPIN(n)$.
}

\subheading{Lemma 2.7}
{\it

When define the degree of $\CLIF(V)$ by a number of elements of $V$,
$\deg(e_1,\cdots,e_n)=n$ and graded Lie bracket,
$$
   [A,B]_\pm = A B - (-)^{\deg(A)\cdot\deg(B)}B A,
$$
$\CLIF_+(V):=\cliffq(\Lambda^2 V)$ with the  graded Lie bracket 
$[,]_\pm$ is a Lie algebra, which is isomorphic to the Lie algebra 
$\roman{so}(V)$, and $\SPIN(V)$ is a compact Lie group associated with
the Lie algebra. 
}

\demo{Proof}(p.105 in [BGV]) 
The isomorphism $\tau:\cliffq(\Lambda^2 V) \to \roman{so}(V)$ is 
obtained by letting $a\in\cliffq(\Lambda^2 V)$
acting on $\cliffq(\Lambda^1 V)\approx V$ by the adjoint action:
$\tau(a)\cdot v=[a,v]_\pm\equiv a v+va$. Then we obtain
$\tau(a): \cliffq(\Lambda^1 V) \to\cliffq(\Lambda^1 V) $.
Hence $:\cliffq(\Lambda^2 V)$ is Lie subalgebra of $\roman{gl}(V)$.
Due to the Jacobi identity,
$$
  [[a,v]_\pm,w]_\pm +[[v,a]_\pm,w]_\pm+[[w,v]_\pm,a]_\pm=0,
  \quad \text{ for } v, w \in V,
$$
and $[[w,v]_\pm,a]_\pm=0$, we obtain 
$(\tau(a)\cdot v,w)=-(v,\tau(a)\cdot w)$.
Hence $\tau(a)$ is isomorphism.
From the construction of $\tau$, $\cliffq(\Lambda^2 V)$
is the Lie algebra associated with the spin group $\SPIN(V)$.
\qed \enddemo

For a matrix $A \in \roman{so}(V)$, the Clifford element
is given as,
$$
     \tau^{-1}(A)= \frac{1}{2}
      \sum_{i,j}(A e_i,e_j) \cliffp(e_i)\cliffp(e_j).
$$

\subheading{Proposition 2.8}
{\it

$\roman{dim} V > 1$, there is a homomorphism $\tau$,
$$
   \tau:\SPIN(V) \to \SO(V), \quad
   ( \tau(g) v = g v g^*),
$$
which is 2-fold covering map.
}
\demo{Proof}(p.106-p.107 in [BGV]) 
$\tau$ is locally isomorphic.
Let $g$ be an element of kernel of $\tau$ or $\tau(g)=1$.
Then for all $v \in \CLIF(\Lambda^1 V)\equiv V$, $[g,v]_\pm=0$.
From $\cliffp([v,g]_\pm)=-2\intp(v) \cliffp(g)=0$,
$g$ must be scalar.
From the definition $g\cdot g^*=1$, the kernel of $\tau$ is
$g = \pm 1$.
\qed \enddemo

\subheading{Corollary 2.9}
{\it
$\pm 1 \in \ZZ_2$ is a center of $\SPIN(V)$.
}

\subheading{Definition 2.10 (Spinor $\CC$ Group)}[G p.186]

{\it

$$
\SPINC(V):= \SPIN(V) \times \CC^\times/(\{ 
   (\epsilon, \epsilon), \epsilon= \pm1\}.
$$
Here we will identify the $\pm1$ as a center of $\SPIN(V)$
and $\pm1 \in \CC^\times$. We will denote it $\SPINC(n)$ for
$V=\RR^n$.
}

Then next is obvious.

\subheading{Proposition 2.11}
{\it
$$
   1 \to \CC^\times \to \SPINC(n) \to \SO(n) \to 1,
$$ $$
1 \to \SPIN(n) \to \SPINC(n) \to \CC^\times \to 1.
$$ 
}

Next we will consider the representation of the spin group.

\subheading{Proposition 2.12 } 

{\it \roster
For $n \in \NN$, these relations hold;

\item 
$$
   \CLIF( \RR^{n+2} ) \otimes_{\RR} \CC \approx 
   \CLIF( \RR^{n} )\otimes_{\RR}\END( \CC^{2} ) .
$$

\item
$$
   \split
   \CLIF( \RR^{2n} ) \otimes_{\RR} \CC &\approx
    \END( \CC^{2^n} ).\\
   \CLIF( \RR^{2n+1} ) \otimes_{\RR} \CC &\approx 
   \END( \CC^{2^n} )\oplus
   \END( \CC^{2^n} ).
   \endsplit
$$
\endroster
}
\demo{Proof}{
Let us define a Pauli matrices 
$\sigma_a \in \END( \CC^{2} )$ ($a = 0,1,2,3$),
$$
   \sigma_0:= \pmatrix 1 & 0 \\ 0 & 1 \endpmatrix, \quad
   \sigma_1:= \pmatrix 0 & 1 \\ 1 & 0 \endpmatrix, \quad
   \sigma_2:= \pmatrix 0 & -\sqrt{-1} \\  \sqrt{-1} & 0 \endpmatrix,
    \quad
   \sigma_3:= \pmatrix 1 & 0 \\ 0 & -1 \endpmatrix, \quad
$$
For an orthonormal basis $\{e_j\}_{i=1,\cdots,n+2}$
of $\RR^{n+2}$, we will introduce
a linear map
$$
   \RR^{n+2} \to  \CLIF( \RR^{n} )\otimes_{\RR}\END( \CC^{2} ),
$$
$$
   e_j \mapsto e_j \otimes \sigma_3, \quad
   e_{n+1} \mapsto 1 \otimes \sigma_2, \quad
   e_{n+2} \mapsto 1 \otimes \sigma_1.
$$
We can extend this relation to $
\EMV(\RR^{n+2}) \to  \CLIF( \RR^{n} )\otimes_{\RR}\END( \CC^{2} )$.
Since
$$
   \sigma_1 \sigma_2 = \sqrt{-1} \sigma_3, \quad 
   \sigma_2 \sigma_3 = \sqrt{-1} \sigma_1, \quad 
   \sigma_3 \sigma_1 = \sqrt{-1} \sigma_2,
$$
$$
   e_j \otimes \sigma_3 \cdot e_k \otimes \sigma_3=
   e_j e_k \otimes\sigma_0.
$$
For $a \in  \CLIF( \RR^{n} )$, following relations hold,
$$
   \split
a\cdot e_{n+1} &\mapsto a \otimes (\sigma_3)^{\deg a} \sigma_2, \\
a\cdot e_{n+2} &\mapsto a \otimes (\sigma_3)^{\deg a} \sigma_1, \\
a\cdot &\mapsto a \otimes (\sigma_3)^{\deg a}  ,\\
a\cdot e_{n+1}e_{n+2} &\mapsto a \otimes (\sigma_3)^{\deg a+1} .\\
\endsplit
$$
So (1) is proved.
On (2), we ill note the relation, as a first step,
$$
   \END( \CC^{n} )\otimes\END( \CC^{2} )= \END( \CC^{2n} ).   
$$
For even case, we can prove it using the relation
$$
   \CLIF( \RR^{2} ) \otimes_{\RR} \CC \approx \END( \CC^{2} ),
$$
which is given by Pauli matrices.
For odd case,
$$
   \CLIF( \RR^{1} ) \otimes_{\RR} \CC \approx \END( \CC )\oplus
   \END( \CC)\approx  \CC \oplus\CC.
$$
$\CLIF( \RR^{1} ) \otimes_{\RR} \CC= \SPAN_\CC\{1,e\} \to
\SPAN_\CC\left\{\dfrac{1+e}{2},\dfrac{1-e}{2}\right\}$.
Inductively we can prove it.
\qed}\enddemo

 \subheading{Proposition 2.13 
  (Spinor Representation of Clifford Algebra)} 
{\it
\roster

\item
$\CLIF(\RR^n) \otimes_{\RR} \CC$ as a $\CLIF(V)$-module is expressed as
$$
   \CLIF(\RR^n) \otimes_{\RR} \CC=
    \sum_{i=1}^{2^{[\frac{n}{2}]}}  \Psi_i ,
$$
where $\Psi_i \approx  \Psi$ for $i=1,\cdots,2^{[\frac{n}{2}]}$.

\item if $n$ is even, $\Psi$ is irreducible as a $\CLIF(\RR^n)$-module.

\item if $n$ is odd, $\Psi\approx \Psi_+\oplus \Psi_-$ and
$ \Psi_\pm$ are irreducible.
\endroster
}

\demo{Proof}{
(1) is obvious from the previous proposition.
If we introduce  a chiral operator and involution
$$
   \Gamma := e_1 e_2 \cdots e_n ,
$$
we have the relation
$$
   \Gamma e_j = (-1)^{n-1} e_j \Gamma .
$$
Hence $\Gamma$ is a center of $\CLIF(V)$ if $n$ is odd. 
Thus $\Gamma v = v \Gamma$ for $v \in \CLIF_+(\RR^n)$ for
any $n\in \NN$.
\enddemo

Finally, we will give only results.

\subheading{Proposition 2.14  (Spinor Representation of Spin Group)} 
{\it
\roster

\item
$\CLIF(\RR^n) \otimes_{\RR} \CC$ as a $\SPIN(n)$-module, called as
spinor-representation, is expressed as
$$
   \CLIF(\RR^n) \otimes_{\RR} \CC
   = \sum_{i=1}^{2^{[\frac{n}{2}]}} \Psi_i ,
$$
where $\Psi_i \approx  \Psi_1$ for $i=2,\cdots,2^{[\frac{n}{2}]}$
and $\Psi_i$ is a representation of $\SPIN(V)$.

\item if $n$ is even, $\Psi$ is
is reduced to $\Psi\approx \Psi_+\oplus \Psi_-$ as $\SPIN(V)$-module
and $ \Psi_+\not\approx \Psi_-$

\item if $n$ is odd, $\Psi\approx \Psi_+\oplus \Psi_-$ and
$ \Psi_+\approx \Psi_-$.
\endroster
}

%\newpage

{\centerline{{\bf \S 3. Spinor Sheaf}}
\tvskip

As we wish to construct a theory in the framework of  the sheaf theory in the
differential geometry and $\DMod$-module, we will go on to study the sheaf
theory for the spinor geometry. As far as I know, there is no reference of
the studies of spinor bundle from sheaf theoretical point of view. However
though Atiyah, Bott and Shapiro did not used terminology of sheaf in the 
study of Clifford module [ABS],  their studies can be easily translated to
language of sheaf theory. In this article, we will translate the theory of
spinor bundle in [ABS, BGV, G, Y] to language of sheaf theory.

Let $(M,\CM)$ be a $n$-dimensional 
real analytic algebraized manifold without singularity
and $\CM$ is its structure
sheaf of real analytic functions over $M$ and it
has the Riemannian metric $(\Theta_M,\fg_M)$,
where $\Theta_M$ is a tangent sheaf [I, Mal1,2].

Generally a sheaf $\Cal S$ (of sets) over a topological space $X$ is 
characterized by a triple $(\Cal S, \pi, X)$ because it is defined so
that there is a local homeomorphism $\pi: \Cal S \to X$ [Mal1]. 
However as we did in [I], we will also write it $\Cal S:=(\Cal S, \pi,
X)$ for abbreviation.
Further we will write the set of (local) sections of $\Cal S$ over an
open set $U$ of $X$ as $\Cal S(U)$ or $\Gamma(\Cal S, U)$. Using the 
category equivalence between category of sheaves and category of 
complete presheaves (due to theorem 13.1 in [Mal1]), we will mix them.

We will go on to use the notations and definitions in the part I [I] 
unless we give notions. The sheaves, $\RR_M$, $\RR^{>0}_M$, $\ZZ_M$,
${\ZZ_2}_M$, $\SHom(\cdot,\cdot)$, $\Theta_M$, $\III(\cdot)$, 
$\SGL(n,\CM)$, $\Omega_M$, $\sqrt{\Omega_M}$, $\sqrt{\CM}$ $\DM$ and so
on, are defined in [I]. We also employ the Einstein convention as in [I].

We will set up the group sheaf as a translation of
a principal bundle as a subsheaf of $\SGL(n,\CM)$
using the results of Mallios [Mal1,2].

\subheading{Definition 3.1}
( Principle Sheaves $\SSSO(n,\CM)$  and $\SSO(n,\CM)$)
[Mal1]{\it\roster

\item 
$SO(M)$ ($O(M)$) is a principle bundle of $\SO(n)$ ($\OO(n)$) over 
$n$-dimensional manifold $M$,
$$\CD        \SO(n)  @>>> SO(M) \\@.  @VVV \\ 
                   @.    M . \\
\endCD$$   

\item $\SSSO(n,\CM)$ ($\SSO(n,\CM)$) is a sheaf for
real analytic sections of
$SO(M)$ ($O(M)$).
\endroster}

For an open set $U$ of $M$, there is a local coordinate system $\{ x^i, 
\partial_i\}_{1\le i \le n}$ and $\{e_j\}_{j=1,\cdots,n}\in \Theta_M$ and 
$\{e^j\}_{j=1,\cdots,n}$ $\in \Omega_M^1$ is an orthonormal set satisfying
the relations $<e^i, e_j> $ $= \delta^i_j$. We will fix the notations in this
section.

\subheading{Definition 3.2} (The Levi-Civita Connection) [Mal1]
{\it

The Levi-Civita connection $\LCconeM_\theta
\in\Gamma(U,\SHom_\CM( \Theta_M, \III(\Theta_M)))$ is defined as an 
integrable connection over the Riemannian tangent sheaf $(\Theta_M,\fg_M)$
for an open set $U$ of $M$ and $\theta \in \Theta_M(U)$ if it holds the 
following relations.

\roster
\item
For sections $\theta_1,\theta_2 \in \Theta_M(U)$,
$$
   \LCconeM_{\theta_1} \theta_2- \LCconeM_{\theta_2} 
   \theta_1=[\theta_1,\theta_2].
$$

\item
When we define $\LCconeM\in\Gamma(U,\SHom_{\RR_M} (\Theta_M, 
\Theta_M\otimes\Omega_M^1))$ as $ \LCconeM:= \LCconeM_{e_j} \otimes e^j$
using the local orthonormal set, 
$$
   d \fg_M(\theta_1,\theta_2) 
   = \fg_M(\LCconeM \theta_1,\theta_2) +\fg_M(\theta_1,\LCconeM  \theta_2).
$$
\endroster
}

\subheading{Lemma 3.3}

{\it
The tangent bundle $TM$ of $M$ can be regarded as an associate bundle for
a representation of $ \rho:\OO(n)\to \RR^n$,
$$
   TM = O(M) \times_{ \roman{O}(n)} \RR^n.
$$
}
\demo{Proof} it is obvious [BGV].\qed \enddemo

\subheading{Corollary 3.4}

{\it 
There is an action of  $\SSSO(n,\CM)$ ($\SSO(n,\CM)$) on $\Theta_M$:
$\SSSO(n,\CM)$ ($\SSO(n,\CM)$) is a subset of $\III(\Theta_M)$.
}

When the Levi-Civita  connection $\LCconeM_\theta$ in $\Gamma(U,\SHom_\CM(
\Theta_M, \III(\Theta_M)))$ belongs to $\Gamma(U,\SHom_\CM( 
\Theta_M, \SSSO(n,\CM)))$ ($\Gamma(U,\SHom_\CM( \Theta_M, \SSO(n,\CM)))$).
Then we call it the  Levi-Civita  connection of $\SSSO(n,\CM)$
 ($\SSO(n,\CM)$).

We also use the correspondence of lemma 3.3 and
 will introduce the Clifford sheaf.

\subheading{Definition  3.5 (Clifford Module)}[ABS, Y]
{\it\roster
 
\item We will define a Clifford bundle $\CLIFs(M)$ as a frame bundle of 
principle $O(M)$ bundle induced from the representation $\roman{O}(n) \to 
\CLIF( \RR^n)$.
$$
   \CLIFs(M) = O(M) \times_{\roman{O}(n)} \CLIF( \RR^n ).
$$
For $p \in M$, its fiber is $ \CLIF(T_p^* M)$.

\item
We will define a sheaf $\Clifs_M$, which
is generated by a set of real analytic sections of $\CLIFs(M)$ over 
an open set $U$ of $M$.
\endroster
}

\subheading{Proposition  3.6 (Clifford Module)}
{\it\roster

\item
When we will regard $\Clifs_M$ as a Clifford module, we have 
decomposition for an open set $U \in M$,
$$
   \Clifs_M(U)=\Clifs_{M+}(U) \oplus \Clifs_{M-}(U),
$$
corresponding  to $\CLIF_+( \RR^n )$ and $\CLIF_-( \RR^n )$
 respectively.
Similarly we have the subbundles of $\CLIFs(M)$ given by the Whitney sum,
$$
  \CLIFs(M)=\CLIFs_+(M)\oplus \CLIFs_-(M).
$$

 \item A sheaf morphism $\cliffq: \Omega_M \to \Clifs_M$ can defined by
 local relation $\cliffq: \Omega_M(U) \to \Clifs_M(U)$.
 
\endroster
}

\demo{Proof}[ABS, Y] Since they are locally defined, they are obvious form
argument of \S2.\qed \enddemo

\tvskip

\subheading{Definition 3.7 (Clifford Connection) }[ABS p.117, Y]
{\it

 An integrable connection $\CLconeM_{ \theta} \in
\Gamma(U,\Hom_\CM( \Theta_M, \III(\Clifs_M))$ for an open set $U \subset M$
and $\theta \in \Theta_M(U)$ is called Clifford connection
if it is satisfied with the following relation:
for $\cliffq(a) \in \Clifs_M(U)$
$$
   [\CLconeM,\cliffq(a)]=\cliffq(\LCconeM a),
$$
where $\LCconeM$ is a Levi-Civita connection.
}

\subheading{Lemma 3.8 (Spin Group Sheaf $\Spin(n,\CM)$) }

{\it
Let us define a  $\SPINs(M)$ as a subbundle $\CLIFs_+(M)$,
 which can be locally defined by the map $\tau : \SPINs(M)|_U \to 
 SO(M)|_U$ of proposition 2.8.
Similarly we will define a  $\Spin(n,\CM)$ as
a group sheaf generated by analytic sections of $\SPINs(M)$,
which is called spin group sheaf.
}

\subheading{Proposition 3.9 (Bockstein exact sequence)}
{\it

For the exact sequence of sheaf,
$$
1_M \to {\ZZ_2}_M\to \Spin(n,\CM)
   \to  \Cal O(n,\CM) \to 1_M,
$$
there is the Bockstein exact sequence of \v Cech cohomology,
$$
   0 \to \check \HH^1({\ZZ_2}_M) \to \check \HH^1( \Spin(n,\CM))
   \to \check \HH^1( \Cal O(n,\CM)) \to \check \HH^2({\ZZ_2}_M).
$$
Here $\check \HH^2({\ZZ_2}_M)$ is called second Stiefel-Whitney class.
}

\demo{Proof} It is obvious from the general theory of 
exact sequence of sheaf. \qed \enddemo

\subheading{Definition 3.10 (Spin $\CC$ Group Sheaf $\SpinC$)}

{\it
The spin $\CC$ group sheaf $\SpinC$ is defined by local sections 
over an open set $U \subset M$
$$
\Gamma(U, \SpinC(n,\CM)):= \Gamma(U, \Spin(n,\CM))
\oplus \CC_M^{\times}(U) /(\{ 
   (\epsilon, \epsilon), \epsilon= \pm1\}).
$$
where $\epsilon $ is a local section of ${\ZZ_2}_M(U)$.
}

\subheading{Proposition 3.11}
{\it

There are  sheaf short exact sequences of group sheaves,
$$
   1_M \to{\CC_M^\times} \to \SpinC(n,\CM) 
   \to \Cal S\SSO(n,\CM) \to 1_M ,
$$
$$
    0 \to \ZZ_M \to \CC_M \overset \exp \to \to \CC_M^\times \to 0.
$$
}

\demo{Proof} From  2.11, there is an exact sequence of their germs
and first one is obvious. Second  is following.
For germ $g_x \in \CC_M^\times$, there exists
germ $f_x \in \CC_M$ such that $g_x =\exp(f_x)$ because
$f_x \equiv  \log g_x $ modulo $2 \pi \sqrt{-1}\ZZ$.
\qed \enddemo

\subheading{Definition 3.12} (Spin $\CC$ structure)
{\it 

We say that $n$-dimensional manifold $M$ has
$\SpinC$ structure if there is its associated principle bundle,
$$\CD       \SPINC(\RR^n)  @>>> \SPINs^{\CC}(M) \\@.  @VVV \\ 
                   @.    M . \\
\endCD
$$
and the cotangent bundle $T^*M$ of $M$  is
equivalent with $ \SPINs^{\CC}(M) \times_{ \SPIN^{\CC}(\RR^n)} \RR^n$.
}

\subheading{Proposition 3.13} (Spin $\CC$ structure) [Y]
{\it

There is a $\SpinC$ structure for a manifold $M$,
if  and only if $\alpha \in \check \HH^2( \ZZ_M)$ such that
$$
   w_2(M) = \alpha \quad \text{\ modulo\ \ }  2\ZZ_M .
$$

} 

\demo{Proof}
Let us introduce a correspondence of germs
$$
   {\tau_1}_x : \SpinC(n,\CM)_x \ni
    [\sigma,\lambda] \mapsto \tau(\sigma) \in \SSSO(n,\CM)_x,
   \quad
    {\tau_2}_x : \SpinC(n,\CM)_x \ni
     [\sigma,\lambda] \mapsto \lambda \in {\CC_M^\times}_x,
$$

First we consider the necessary condition. Suppose that
$M$ has $\SpinC$ structure. By $\tau_2$, we regard ${\CC_M}^\times$ as 
the $\SpinC$ module and as the determinant module of $\SpinC$.
By the Weil-Kostant theorem [Br p.66], we call 
$c \in \check \HH^1(\CC_M^\times)$
the first Chern class and it is isomorphic to the Picard group 
$\check \HH^2(\ZZ_M)$. 

For a local section $[\sigma_i, \lambda_i]$ of $\SpinC(n,\CM)$ over $U_i$,
we will express the local map over $U_{ij}\equiv U_i\cap U_j\neq \empty$,
$
   [\sigma_i, \lambda_i]=[\sigma_{ij},
    \lambda_{ij}][\sigma_i, \lambda_i]
$
where $[\sigma_{ij}, \lambda_{ij}]\equiv
=[\sigma_i\sigma_j^{-1}, \lambda_i/\lambda_j]
$.
Then for $U_{ijk}\equiv U_i\cap U_j\cap U_k \neq\empty$,
we have the relation from the definition of $\SpinC$,
$$
       \sigma_{ij}\sigma_{jk}\sigma_{ki}=\lambda_{ij}
       \lambda_{jk}\lambda_{ki}.
       \tag *
$$
Due to the relation $\tau(\sigma_{ij}\sigma_{jk}\sigma_{ki}) =1$, 
$\sigma_{ij}\sigma_{jk}\sigma_{ki}$ must be an element of $\ZZ_2=\{\pm
1\}$ valued second \v Cech cohomology, {\it i.e.}, element of 
Stiefel-Whitney class $w_2(M)$. Let $\epsilon_{ijk} $, 
$[\epsilon_{ijk}]\equiv w_2(M)$ be the value of above relation (*). The
left hand side must be an element of the Picard group. We have the 
relation $w_2(M) = \alpha$ modulo $2\ZZ_M$. 

The sufficiency is given as follows.
We take the element $\alpha$ of the Picard group $\check \HH^2(\ZZ_M)$
holding the relation in the propositions for all open set of $M$. Then
we have a line bundle $LINE(M)$, whose first Chern class is $\alpha$,
so that for an open set $U$ of $M$, its local trivialization is consist
with $T^*M|_U= (\SPINs^{\CC}(M) \times_{ \SPIN^{\CC}(\RR^n)} \RR^n)|_U$
by choosing the relation (*). We can find the set 
$\{[\sigma_{ij},\lambda]\}$ consisting with the relation (*) for all 
open set $U_{ijk}$. Then we have the spin $\CC$ structure. \qed \enddemo

\subheading{Definition 3.14 (Spin Representation)} [G, Y]

{\it
$M$ is a real $n$-dimensional compact manifold $M$
which has $\SpinC$ structure. 
Let an open set $U \subset M$,  
$$
   \Psi(U) = \SPINCs(U) \times_{\SPINC(\RR^n)} \Psi(\RR^n).
$$

Let $\varPsi_M$ be a sheaf of analytic section
of $\Psi(U)$ of $U \subset M$, which is called spinor sheaf.
}

Then we have a relation $\varPsi_M$ is sheaf isomorphic to 
$\CC_M\otimes_{\RR_M}\Clifs_M$.

\subheading{Definition 3.15 (Levi-Civita connection of $\SpinC$)}[G, Y]

{\it
Suppose that a real $n$-dimensional compact ringed manifold
$(M,\CM)$  has the Riemannian structure $(\Theta_M,\fg_M)$ 
and $\SpinC$ structure. 
If the connection of $\SpinC$ is induced from the
the Levi-Civita connection $\LCconeM$ of $SO(n, \CM)$, we will call
it the Levi-Civita connection of $\SpinC$.
}

\subheading{Proposition 3.16 (Levi-Civita connection of $\SpinC$)}[Y]
{\it

The Levi-Civita connection $\SPconeM$ of $\SpinC$  is 
the Clifford connection $\CLconeM$.
}

\demo{Proof} From the corollary 3.4 and the correspondence of $\tau_1$
in the proof of proposition 3.13, it is obvious. \qed \enddemo

\subheading{Definition 3.17 (Dirac Operator)}[G, Y]

{\it
The Dirac operator $\Nabla_M^{(\Dir)} \in \III(\Psi) $
 is defined by
$$
   \Nabla_M^{\Dir}
    = \sum_{i=1}^n \cliffq(e^i) \cdot \SPconeM_{\theta_i} ,
$$
where $<e^i, \theta_j > = \delta^i_{\ j}$.
If we denote $\SPconeM_{\partial_i} (e^j) = \sum_k w^{j}_{i k } e^k$,
$$
    \SPconeM_{\partial_i} = \partial_i +
    \frac{1}{4} \sum_{j k }w^{j}_{i k } \cliffq(e^k)\cliffq(e^j) .
$$
}

\subheading{Corollary 3.18 (Dirac Operator in $\EE^n$)}[G, Y]

{\it
The Dirac operator in $\EE^n$ is given by
$$
   \Nabla_{\EE^n}^{\Dir} = \sum_{i=1}^n \cliffq(d x^i)
    \cdot \partial_i ,
$$
where $x^i$ is the Cartesian coordinate. 
}

We sometimes use the notations  $\Partial:=\Nabla_M^{\Dir} $ [Mat11, 16].

\subheading{Definition 3.19 (Dirac System)}[I-Bjo]

{\it

We introduce a coherent $\DM$-module for
the Dirac equation, which is locally expressed by an
exact sequence
$$
   \DM^{2^{[\frac{n}{2}]}} 
   \overset \sqrt{-1}\Nabla^{\Dir}\to{\longrightarrow}
    \DM^{2^{[\frac{n}{2}]}}  \overset \eta\to{\longrightarrow} 
   \DDM_M^{\Dir}\to 0.
$$
Let this $\DM$-module $\DDM_M^{\Dir}$ be referred Dirac system in this 
article. }

\tvskip
Here for $\eta( \epsilon_i) = u_i$ for
base $\epsilon_i := ( 0, \cdots, 0, 1, 0,\cdots 0)\in 
\DM^{2^{[\frac{n}{2}]}}$, $\DDM_M^{\Dir}$ is expressed as 
$\DDM_M^{\Dir}= \sum \DM u_i$. As the Dirac operator is acted by group
action of $\SpinC$, $\SpinC$ acts on $\DDM_M^{\Dir}$.

As we argued in \S 3 in [I], we will also consider the 
(anti-)self-adjoint
operators for the spinor representations module $\varPsi_M$.
We have the Hodge-star operator $*$ and it acts on $\varPsi_M$
$$
   *:\varPsi_M \to \varPsi_M \otimes_{\CM} w_M ,
$$
where $w_M$ is the volume form of $M$.

Let the local coordinate system $\{x^i\}_{1\leq i \le n}$
have non-trivial Riemannian metric section 
$\fg_M = g_{M i,j} d x^i \otimes d x^j$,
$g_M:= \det g_{i,j}$ and $w_M=g_M dx^1 \cdots dx^n$.
Then in general $\partial_i$ is not self-adjoint for $ \sqrt{\CM}$
because $*:\sqrt{\CM^\CC} \to \sqrt{\CM^\CC}w_M$ and adjoint
form of $\partial_i$ also acts upon $w_M$. Thus we
can define the anti-self-adjoint operator  ${\tnablaM}_i
= g_M^{-1/4} \partial_i g_M^{1/4}$.
Further we introduced a half form $\sqrt{\omega_M}$ in [I].

\subheading{Definition 3.20 (Spin Half Form \& Spin Wave Sheaf)}

{\it\roster 

\item 
Let $\varPsi_M\otimes_{\RR_M} \sqrt{\omega_M}$ be spin half form.
The local section of 
connection of spin half form 
$\SHom_M( \Theta_M,\III(\varPsi_M\otimes_{\RR_M}\sqrt{\omega_M}))$ 
is expressed by ${\SPHcone_M}_\theta$ for $\theta\in \Theta_M(U)$:
$$
    {\SPHcone_M}_{\partial_i} = \partial_i 
    - \frac{1}{4} \partial \log g_M
    +\frac{1}{4} \sum_{j k }w^{j}_{i k } \cliffq(e^k)\cliffq(e^j) .
$$
   
\item 
Let $\varPsi_M\otimes_{\RR_M} \sqrt{\CM^\CC}$ be Spin wave sheaf.
The local section of 
connection of spin half form $\SHom_M( \Theta_M,
\III(\varPsi_M\otimes_{\RR_M}\sqrt{\CM}))$ 
is expressed by ${\SPWcone_M}_\theta$ for $\theta\in \Theta_M(U)$:
$$
    {\SPWcone_M}_{\partial_i} = {\tnablaM}_i - \frac{1}{4} 
    \partial \log g_M
    +\frac{1}{4} \sum_{j k }w^{j}_{i k } \cliffq(e^k)\cliffq(e^j) .
$$
where    ${\tnablaM}_{\partial_i}$ is 
locally defined by ${\tnablaM}_i= g_M^{-1/4} \partial_i g_M^{1/4}$
using the local coordinate system.
\endroster
}

\subheading{Proposition 3.21 (Spin Half Form)}

{\it
The category whose object and morphism are given by
$(\varPsi_M\otimes_{\RR_M} \sqrt{\omega_M}, {\SPHcone_M}_\theta)$
is category equivalent with that of 
 $(\varPsi_M\otimes_{\RR_M} \sqrt{\CM},{\SPWcone_M}_\theta)$.
}

\demo{Proof} This is essentially same as the proposition 2.5 in [I].
\qed \enddemo

\subheading{Remark 3.22}

In fashion of [I], we should discriminate  ${\SPWcone_M}_{\partial_i}$
and ${\SPconeM}_{\partial_i}$ and  explicitly write action of 
$\overline{\sigma}$ of the correspondence
defined in the lemma I-3.24.
However since 
${\SPWcone_M}_{\partial_i}$ $\equiv {\SPconeM}_{\partial_i}$ in
the primary meanings,
we will neglect $\overline{\sigma}_M$ in this article for
abbreviation.

%\newpage

{\centerline{\bf{\S 4.  Dirac Operators 
in a submanifold $S \hookrightarrow\EE^n$}}
\fvskip

Let $S$ be a $k$-dimensional real  analytic compact submanifold of 
$\EE^n$, associated with an integrable connection $\tnablaS^{(s)}$, 
embedded in $\EE^n$;
$\iota_S: S \hookrightarrow \EE^n$.
Let $T_S$ a tubular neighborhood of $S$, $i_S :S \hookrightarrow T_S$
and $ i_{T_S}: T_S \hookrightarrow \EE^n$ such that
 $\iota_S \equiv i_{T_S} \circ i_S$.
$T_S$ has a projection $\pi_{T_S} : T_S \to S$.
$\EE^n$ has a natural metric $\fg_{\EE^n}$ and thus $T_S$ and $S$ has
Riemannian module $(i_{T_S}^{-1} \Theta_{\EE^n}, i_{T_S}^{-1} \fg_{\EE^n})$
and $(\Omega_S, \iota_{S}^{-1} \fg_{\EE^n})$.

\fvskip

\subheading{Notations 4.1}
{\it 

\roster

\item The  inverse image of the Dirac system of  ${\DDM}_{T_S}^{\Dir}$ 
is defined as
$$
   {\DDM}_{S \to T_S}^{\Dir}= 
   \DS \otimes_{i_S^{-1} \DT} i_S^{-1}({\DDM}_{T_S}^{\Dir}),
$$
where
$$
  {\DDM}_{T_S}^{\Dir}=\DT \otimes_{i_{T_S}^{-1} \DMod_{\EE^n}} 
  i_{T_S}^{-1}({\DDM}_{\EE^n}^{\Dir}).
$$

\item The tangent sheaf of $\Theta_{T_S}$ is $\Theta_{T_S}=i_{T_S}^{-1} 
\Theta_{\EE^n}$ and has a direct decomposition as a $\CTS$-module,
$$
   \Theta_{T_S} = \Theta_{T_S}^\parallel\oplus \Theta_{T_S}^\perp,
$$
where $i_S^{-1}\Theta_{T_S}^\parallel =\Theta_{S}^\parallel$, $
i_S^{-1}\Theta_{T_S}^\perp =\Theta_{S}^\perp$
and $\Theta_{S}^\perp$ is
the normal sheaf defined by the exact sequence,
$$
 0 \to \Theta_S \to \iota_S^{-1} \Theta_M \to  \Theta_{S}^\perp \to 0.
$$

\item An affine vector (coordinate) $\bold Y\equiv(Y^i)$ in $T_S \subset 
\EE^n$  is expressed by,
$$
    \bold Y= \bold X + \bee_\dalpha q^\dalpha,
$$
for a certain affine vector $\bold X$ of $S$. 

\item A point $p$ in $T_S$  is expressed by the  local coordinate 
$(u^\mu):=(s^1,s^2,$ $\cdots,s^k,$ $q^{k+1},\cdots,$ $q^n)$, 
$\mu=1,2,\cdots,n$ where $(s^1, \cdots,s^k)$ is a local coordinate of 
$\pi_{T_S} p$; We assume that the beginning of the Greek ($ \alpha$, $\beta$,
$\gamma$, $\cdots$) runs from 1 to $k$ and they with dot ($ \dalpha$, 
$\dbeta$, $\dgamma$, $\cdots$) runs form $k+1$ to $n$). Let 
$(u^\mu)=(s^\alpha,q^\dalpha)$, where the ending of the Greek ($ \mu$, $\nu$,
$\lambda$, $\cdots$) run from $1$ to $n$.

\item $\bEE_\mu:=\partial_{\mu}:=\partial/ \partial u^{\mu}$ 
is a base of $ \Theta_{T_S}(U)$.  
$\bEE^\mu:=d u^{\mu}$ 
is a base of $ \Omega^1_{T_S}(U)$: 
$<\bEE^\mu, \bEE_\nu>= \delta^\mu_\nu$

\item 
For $U \subset T_S$, the induced metric of $T_S$ from $\EE^n$ has a direct
sum form,
$$
   \fg_{T_S}:=i_{T_S}^{-1}\fg_{\EE^n}=\fg_{T_S^\parallel}
   \oplus\fg_{T_S^\perp},
$$
where $\fg_{T_S^\perp}$ has a trivial structure.
In local coordinate, 
$$
\fg_{T_S^\perp}= \delta_{\dalpha,\dbeta} \dd q^\dalpha\otimes
 \dd q^\dbeta, \quad 
\fg_{T_S^\parallel}
     =g_{T_S\alpha\beta} \dd s^\alpha \otimes \dd s^\beta ,
$$
or for $g_{T_S\mu,\nu}:=\fg_{T_S}(\partial_\mu,\partial_\nu)$,
$$
 g_{T_S\dalpha\dbeta}=\delta_{\dalpha\dbeta},
     \quad
     g_{T_S\dalpha\beta}=g_{T_S\alpha\dbeta}=0,
$$
where $\partial_\mu:=\partial/\partial u^\mu$.

\roster }

\fvskip

\subheading{Proposition 4.2}
{\it

 Let the anti-self-adjoint connection 
$\tnablaS_\dalpha\in\Gamma(U,i_S^*\XiT(\Theta_{T_S}^{\perp}))$ 
 for  an open set $U$ in $S$.
There is an injective endomorphism of 
the Dirac system ${\DDM}_{S \to T_S} $,
$$
   \eta^\conf_\dalpha 
     :{\DDM}_{S \to T_S}^{\Dir} \to {\DDM}_{S \to T_S}^{\Dir} ,
$$
for $P \in \Gamma(U,{\DDM}_{S \to T_S}^{\Dir} ) $,
$
       \eta^\conf_\dalpha(P) =  P \tnablaS_\dalpha \in
        \Gamma(U,{\DDM}_{S \to T_S}^{\Dir} ).
$
Then we have a submodule of ${\DDM}_{S \to T_S}^{\Dir} $, 
$$
\eta^\conf:(\overline{\SS}_{S \to T_S} )^{n-k}
\to \sum_{\dalpha=k+1}^n{\DDM}_{S \to T_S}^{\Dir} 
\tnablaS_\dalpha \subset {\DDM}_{S \to T_S}^{\Dir} .
$$

}

\demo{Proof} For the tubular neighborhood, we have a local coordinate
system as we demonstrated in [I]. Using such local coordinate system, 
this can be proved as we did in proposition I-4.1. \qed \enddemo

\subheading{Definition 4.3}
{\it
\roster

\item
We will define a coherent $\DMod_{S}$-module,
${\DDM}_{S\hookrightarrow \EE^n}^{\Dir}$, by the exact sequence,
$$
   ({\DDM}_{S \to T_S}^{\Dir})^{n-k} \overset{\sigma^\conf
   }\to{\longrightarrow}{\DDM}_{S\to T_S}^{\Dir} 
   \longrightarrow
   {\DDM}_{S\hookrightarrow \EE^n}^{\Dir}\longrightarrow 0 ,
$$
Let us call it  submanifold Dirac system 
of $S\hookrightarrow \EE^n$.

\item 
Let the submanifold Dirac system 
${\DDM}_{S\hookrightarrow \EE^n}^{\Dir}$
be decomposed by the exact sequence,
$$
\DMod_S^{2^{[\frac{n}{2}]}} \overset \sqrt{-1}
 {\Nabla}_{S\hookrightarrow \EE^n}
\to{\longrightarrow} \DMod_S^{2^{[\frac{n}{2}]}} \longrightarrow 
   {\DDM}_{S\hookrightarrow \EE^n}^{\Dir}\to 0.
$$
where we tune ${\Nabla}_{S\hookrightarrow \EE^n}$ such that
there is a map $\sigma: \Clifs_S \to i_{S}^{-1} \Clifs_{T_S}$ 
$$
        [\sigma({\Nabla}_{S}) - {\Nabla}_{S\hookrightarrow \EE^n} ] \in 
        \CS \otimes i_{S}^{-1}\Clifs_{T_S}.
$$
we will call
${\Nabla}_{S\hookrightarrow \EE^n}$ is the
 the submanifold Dirac operator of $S\hookrightarrow \EE^n$.

\endroster
}

We will describe our main theorem, which is given in 
[Mat1-3, 5, 6, 10, 11, 16].

\subheading{Theorem 4.4}

{\it
\roster

\item $k=1$ and $n=3$ case, $S$ is curve $C$,
$$
   \Nabla_{C \hookrightarrow \EE^2}^{\Dir} 
   = \pmatrix \partial_s & \frac{1}{2} \kappa_\CC \\
    -\frac{1}{2} \overline{\kappa_\CC} & \partial_s \endpmatrix ,
$$
where $s$ is the arclength of the curve $C$,
 $\kappa_\CC$ is the complex curvature of $C$, defined by
$\kappa_\CC = \kappa(s) \exp\left(\sqrt{-1} 
\dsize\int^s \tau ds\right)$ using 
the Frenet-Serret curvature $\kappa$ and torsion $\tau$.

\item $k=2$ and $n=3$ case, $S$ is a conformal surface,
$$
      \Nabla_{S \hookrightarrow \EE^2}^{\Dir} = 
       2\rho^{1/2}\pmatrix  p & \partial \\
       \overline \partial & p  \endpmatrix\rho^{-1/2} ,
$$                                           
where $(z,\bar z)$ is a complex parameterization of $S$,
$p=\rho^{1/2} H/2 $  
when the volume element of $S$  
is given by dvol $= \rho d z d \bar z$ and 
$H$ is the mean  curvature.

\item $k=2$ and $n=4$ case, $S$ is a conformal surface,

$$
      \Nabla_{S \hookrightarrow \EE^4}^{\Dir} =  2\rho^{1/2}
        \pmatrix & &  p_c & \partial \\
      &  & \overline \partial & -\overline p_c \\
              \overline p_c & \partial & & \\
           \overline \partial &  -p_c & & \endpmatrix
           \rho^{-1/2},
$$
where $(z,\bar z)$ is a complex parameterization of $S$,
$p_c = \rho^{1/2} H_c$  
when the volume element of $S$  is dvol $= \rho d z d \bar z$
and $H_c$ is "complex" mean curvature (see example 4.14 (3) in [I]).

\endroster
}

\subheading{Remark 4.5}

\roster
\item
The $(k,n)=(1,2)$ case is studied by me. When the affine vector of $C$ 
($i_C:C \hookrightarrow \EE^2$) is expressed by $X= X_1 + \sqrt{-1} X_2$,
we have a relation,
$$
       \Hom_\CC(\DDM_{C \hookrightarrow \EE^2}^{\Dir},
       i_C^{-1}C_{\EE^2}^\omega)
   = \CC\pmatrix \sqrt{\partial_s X}\\ 
   \sqrt{-\partial_s X} \endpmatrix.
$$
This is essentially the same as the Frenet-Serret relation.

\item The case that $(k,n)=(2,3)$ and $S$ is a conformal surface is 
identified with the generalized Weierstrass relation [Mat11] whose form was
found by Konopelchenko and Taimanov [Ko1, 2, KT, T1, 2]. (Kenmotsu found more
naive expression of the generalized Weierstrass relation [Ke]). Let $\EE^3
\approx \CC \times \EE\ni (Z:= X_1 + \sqrt{-1} X_2,X^3)$, the solution of the
Dirac equation is given as,
$$
       \Hom_\CC(\DDM_{S \hookrightarrow \EE^3}^{\Dir},
       i_S^{-1}C_{\EE^3}^\omega)
   = \CC\pmatrix \sqrt{\partial_s Z}\\ 
   \overline{\sqrt{-\partial_s Z}}\endpmatrix.
$$

\item 
Konopelchenko found the generalized Weierstrass relation as the zero 
mode of the Dirac operator $\DDM_{S \hookrightarrow \EE^3}^{\Dir}$ 
through the studies of geometrical interpretation of higher dimensional
soliton. Bobenko also pointed out that the Dirac operator plays an 
important role in geometrical interpretation of soliton equation [Bo].
It known that Abresh already found similar equation ending of 1980's on
the investigation of constant mean curvature surface. 

In [Mat2], I conjectured that the submanifold Dirac operator might exhibit
submanifold feature and be connected to higher dimensional soliton
as open problems. 
Berguress and Jenssen followed my approach [TM] and calculated the Dirac
operator of a surface immersed in $\EE^3$ [BJ]. However their surface
is not conformal, it was too complex for me to obtain any meaningful 
result. As I read [KT], I got the above result [Mat11] by means of more
physical method.

On the other hand, Pedit and Pinkall developed the quaternion algebraic
geometry and  reformulated the Dirac operator on a conformal surface in
$\EE^3$ in  framework of quaternion algebraic geometry [Pe,PP]. 
Friedrich also gave a the Dirac operator of a surface immersed in 
$\EE^3$ through the study of spinor bundle [Fr]. He gave a spinor 
connection of a submanifold and his approach is very resemble 
to our methods.
However it is not clear why one considers the Dirac operator since his
approach can not connect with a submanifold Schr\"odinger operator. We
must emphasize that  submanifold differential operators are universal
and  not a special feature of  the spinor bundle. In fact, in optics or
sound wave, we encounter similar differential effects [Mi]. Accordingly
as we described in the \S 1 introduction,the Dirac operator $\DDM_{S 
\hookrightarrow \EE^n}^{\Dir}$ should be studied with the Schr\"odinger
operator $\Delta_{S \hookrightarrow \EE^n}^{\Dir}$ and so on.
 
Further I also gave more concrete conjecture on the
Dirac operator on a conformal surface in $\EE^4$ in ending of 
1997:
the Dirac operator $\Nabla_{S \hookrightarrow \EE^4}^{\Dir}$
is the operator of the  generalized Weierstrass relation of 
a conformal surface in $\EE^4$ [Mat16].  
At the same time, Konopelchenko [Ko2]
and Pedit and Pinkall [PP]  
investigated the  generalized Weierstrass relation of 
a conformal surface in $\EE^4$ and found the 
Dirac operator $\Nabla_{S \hookrightarrow \EE^4}^{\Dir}$.
In other words, as things turned out, my conjecture was proved.

\endroster

Here after we will prove  the theorem 4.4.

\subheading{Lemma 4.6}

{\it
\roster

\item
A element of left-hand $\DMod_S$-module 
$P\in {\DDM}_{S\to T_S}^{\Dir}$ is
also given by $Q\in{\DDM}_{T_S}^{\Dir}$
such that $P := Q|_{ q^\dalpha\equiv 0, \dalpha = k+1,\cdots,n)})$.

\item
An element of left-hand $\DMod_{S}$-module 
$P\in {\DDM}_{S\hookrightarrow \EE^n}^{\Dir}$ is given by
$Q\in{\DDM}_{S \to T_S}^{\Dir}$ such that
$P:=Q|_{\tnablaS_\dalpha\equiv0, (\dalpha = k+1,\cdots,n)}$.

\item The Dirac operator  $\Nabla_{S \hookrightarrow \EE^n}^{\Dir}$
is unique up to the action of $\SpinC$.

\endroster
}

\demo{Proof} As did in proposition I-4.6, we can prove them. \qed\enddemo 

\subheading{Lemma 4.7}

{\it
\roster

\item
The Dirac  operator 
$\Nabla_{T_S}^{\Dir} \equiv i_{T_S}^{-1}\Nabla_{\EE^n}^{\Dir} $, 
which  is expressed by affine coordinate 
$ \Nabla_{T_S}^{\Dir} = \cliffq(d X^i) \cdot \SPcone_{T_S \partial_i} $,
can be written by a local coordinate system,
$$
   \Nabla_{T_S}^{\Dir} =  \cliffq(E^\mu)\cdot {\SPcone_{T_S}}_\mu ,
$$
where  

$$
   {\SPcone_{T_S}}_\mu\equiv {\SPcone_{T_S}}_{E_\mu}
    =\partial_\mu + {\Omega_{T_S}}_\mu,
$$
$\partial_\mu\equiv E_\mu$, $d u^\mu \equiv E^\mu$,
$
{\Omega_{T_S}}_\mu=\frac{1}{4}{\Omega_{T_S}}_{\mu\nu}^\lambda 
\cliffq(E^\mu)\cliffq(E_\lambda)$
and ${\LCcone_{T_S}}_\mu E_\nu 
={\Omega_{T_S}}_{\mu\nu}^\lambda E_\lambda.$

\endroster}

\demo{Proof} Direct computations leads the above results. \qed \enddemo

\subheading{Lemma 4.8}

{\it

\roster

\item
$$
   {\SPcone_{T_S}}_\alpha = \pi_{T_S}^*{\tnablaS}_\alpha+ 
    {\Omega_{T_S}}_\alpha-\frac{1}{4} 
   \partial_\alpha \log g_{T_S}.
$$

\item
$$
   {\SPcone_{T_S}}_\dalpha = \tnablaT_{ E_\dalpha}-\frac{1}{4} 
   \partial_\dalpha \log g_{T_S}.
$$

\item For $k=2$ case,
$$
-\frac{1}{4}\partial_\dalpha \log g_{T_S} =
     \frac{ H_{\dalpha-2} -\partial_\dalpha J(q^\dbeta)/2 }
   {1-2 H_\dalpha q^\dalpha+J(q^\dbeta)}. 
$$

\endroster
}

\demo{Proof} Using the expressions in propositions 
and lemmas 4.5-10 in [I], above results are 
obtained [BJ, Mat11, 16]. \qed \enddemo

\subheading{Lemma 4.9}

{\it

\roster
\item When we define $\Nabla_{S \to T_S}^\Dir
:=\Nabla_{T_S^\parallel}^\Dir|_{ q^\dalpha\equiv0}$, 
$\Nabla_{S\hookrightarrow \EE^n}^\Dir$ is expressed as
$
\Nabla_{S\hookrightarrow \EE^n}^\Dir$  $=
\Nabla_{S \to T_S}^\Dir|_{\tnablaT_{ E_\dalpha}\equiv 0}$.

\item For $k=2$ case, 
$\Nabla_{S\hookrightarrow \EE^n}^\Dir$ is expressed by 
$
\Nabla_{S\hookrightarrow \EE^n}^\Dir$ $=\Nabla_{S}^\Dir 
+\cliffq( E^\dalpha) H_\dalpha
$.
\endroster
}

\demo{Proof} From the definition 4.3, (1) is obvious.
 Using (1) and lemma 4.7 (3),
(2) is obtained.
\qed \enddemo

For a while, we will deal with only $n=4$ and $k=2$ case.

\subheading{Lemma 4.10} [Po, Mat16]

{\it
\roster

Let $(n,k)=(4,2)$ and $S$ be a conformal surface and $\rho:= 
g_{S}^{1/2}$. A natural euclidean inner space coordinate system is 
expressed by $y^a,y^b,\cdots$, $(a,b=1,2)$. The moving frames of $S$ is
given by $d s^\alpha \equiv e^\alpha$, $\partial_\alpha \equiv e_\alpha$.
 
\item
The moving frame is written as,
$$
  e^a_{\ \alpha}=\rho^{1/2} \delta^a_{\ \alpha}.
$$

\item
The Christoffel symbol over a conformal surface $S$ is calculated
as,
$$
   \gamma^\alpha_{\ \beta \gamma}=
   \frac{1}{2} \rho^{-1}
     (\delta^\alpha_{\ \beta}\partial_\gamma \rho
      +\delta^\alpha_{\ \gamma}\partial_\beta \rho
      -\delta_{\beta \gamma}\partial_\alpha \rho) .
$$

\item
The spin connection over $S$ 
$\omega_{\alpha}:=\Omega_{\alpha}|_{q^\dalpha\equiv0}$ becomes
$$
   \omega_{\alpha}=- \frac{1}{4} \rho^{-1}\sigma^{ab}
   (\partial_a  \rho \delta_{\alpha b}
   -\partial_b  \rho \delta_{\alpha a}) , 
$$
where $\sigma^{ab}:=1\otimes [\sigma^a,\sigma^b]/2$. 

\item
The
Dirac operator $\Nabla_{S}^{\Dir}$ can be expressed as
$$
\Nabla_{S}^{\Dir}
=\sigma^1\otimes\sigma^a \delta^\alpha_{\ a}
   [\rho^{-1/2}\partial_\alpha 
   + \frac{1}{2}\rho^{-3/2} 
   (\partial_\alpha \rho)]  .
$$

\item 
The anti-self-adjoint operator of $S$ is given by 
$$
   \tnablaS_{ \alpha}= \rho^{-1/2} \partial_\alpha \rho.
$$

\item
The normal Clifford elements can be connected to the Pauli matrices
$\sigma^a$,
$$
    \cliffq(e^\alpha)=e_a^{\ \alpha} \sigma^
    1 \otimes_\CC  \sigma^a, \quad
    \cliffq(e^\dalpha)= \sigma^2\otimes_\CC \sigma^{\dalpha-2}.
$$

\endroster
}

\demo{Proof} See the p.232 in [Po] or direct computations also give 
(1)-(5). (6) was obtained in [Mat11,16] by direct computations along the
definition of 4.3 \qed \enddemo

\subheading{Corollary 4.11}

{\it

$\Nabla_{S\hookrightarrow \EE^4}^\Dir$ can be expressed by
$$
\Nabla_{S\hookrightarrow \EE^4}^\Dir=
\rho^{-1}[\sigma^1\otimes\sigma^a \delta^\alpha_{\ a}\partial_\alpha 
+ \rho^{1/2}\sigma^1\otimes\sigma^a H_a ]\rho^{1/2}
  .
$$
}

\subheading{Lemma 4.12}

{\it \roster
\item
Let $S\subset\EE_s^3\times \EE_t^1\equiv\EE^4  $ be a product manifold 
 $C \times \EE_t^1$ such that for $\pi_{\EE^4\to \EE_s^3}$ 
 $: \EE^4 \to \EE_s^3$,
$\pi_{\EE^4\to \EE_s^3}(S)$ is a space curve $C$ in $\EE^3$. We set
$s$ be an arclength of $C$ and then we obtain the relation,
$$
   H_1= -\frac{1}{2} \text{tr}_2(\gamma^1_{\ 31}),\quad
   H_2= -\frac{1}{2} \text{tr}_2(\gamma^1_{\ 41}), \quad
   \rho=1. 
$$

\item
Let consider the case that the surface $\SS$  immersed in $\EE^3$ and  the
direction ${\bold e}_4$ is globally identified with the direction of $X^4$.
It means that $H_2= 0$ and $p_c$ is obtained as,
$$
   p_c \equiv p:= \frac{1}{2}\sqrt{\rho} H_1 .
$$
}

\demo{Proof} They are obvious\qed \enddemo

\subheading{Proof of Theorem 4.4}

Corollary 4.11 and lemma 4.12 gives a complete proof of theorem 4.4. \qed

\vskip 0.5 cm

Here we will mention the results in [Mat2, 6, 14].
Though we will not give proofs in this article, following results were 
also obtained by physical computations.

\subheading{Proposition 4.13}

{\it
\roster

\item $S$ is one dimensional case, ${\Nabla}_{S\hookrightarrow \EE^n}^{\Dir}$
is identified with the Ferret-Serret Relation and one of Lax operators of the
generalized MKdV equations .

\item $S$ is a conformal surface  and $n=3,4$ case, 
${\Nabla}_{S\hookrightarrow \EE^n}^{\Dir}$ is identified with the 
generalized Weierstrass relation and related to generalized MNV equation
.
\endroster
}

\demo{Comment} (1) was obtained in [Mat1-3, MT] and (2) was proved by 
another approach in [Ko1,2, KT, T1,T2]. \enddemo

\subheading{Proposition 4.14}
{\it 
\roster

Let us consider a localized ring
$\Cal C_S^{\omega\CC}[[ {\Nabla}_{S\hookrightarrow \EE^n}^{\Dir}]]$ 
and consider
$\det({\Nabla}_{S\hookrightarrow \EE^n}^{\Dir})$.

\item $S$ is one dimensional case, 
$\log \det {\Nabla}_{S\hookrightarrow \EE^n}^{\Dir}$ is 
identified with the
Euler-Bernoulli functional $\Cal B = \int \dd s  |k|^2$.

\item $S$ is one dimensional case, 
$\roman{index} {\Nabla}_{S\hookrightarrow \EE^n}^{\Dir}$ gives the 
fundamental group of the loop .

\item $S$ is a conformal surface case, 
 $\log \det {\Nabla}_{S\hookrightarrow \EE^n}^{\Dir}$ is related to the 
 Willmore functional 
 $\Cal W = \int \dd z \dd \bar z \rho H^2$ and area $\Cal
 A$ of the surface. \endroster
}

\demo{Comment} (1) was obtained in [Mat2, 6, 8], (2) was in [Mat2, 6, 8] and 
(3) was in [Mat14]. (1) and (3) were proved by path integral method and
heat kernel method. (3) was calculated using the path integral and 
generalized Hurwitz $\zeta$-regularzation.  \enddemo

\subheading{Remark 4.15}

For the case of $n=3$ dimensional case, 4.14 (3) means that
$$
   \Cal W= \mu^2 \Cal A +\frac{5}{3}\pi \chi 
   -\nu -\log \det {\Nabla}_{S\hookrightarrow \EE^n}^{\Dir},
$$
where $\mu$ and $\nu$ are certain real numbers and $\chi$ is Euler number.
Hence the Willmore conjecture should be studied in the framework of this
Dirac operator [Mat14].

\subheading{Conjecture 4.16}
{\it 

For any manifold $M$ and $S$, the Dirac operator can represents the
submanifold system. 
} 

\subheading{Remark 4.17}

We should note that differential algebra sheaf $\DMod_S$ generated by
$\tnablaS_\theta$ ($\theta\in \Theta_S$) commutes with that 
$\DMod_S^\perp$ generated by $\tnablaS_\theta$ ($\theta\in 
\Theta_S^\perp$). In the algebra $\DMod_{S\to \EE^n}$, these are 
commutant if we used terminology in von-Neumann algebra. Then $S$ itself
might be regarded as $\DMod_{S\to \EE^n}\cap \Theta_S $ or factor of the 
system. Hence above behavior might be very natural in non-commutative algebra
theory.

\subheading{Remark 4.18}

Since our Dirac operator  represents the submanifold completely,
at least, surface and curve case, we can naturally determine
the operator form of ${\Nabla}_{S\hookrightarrow \EE^n}^{\Dir}$
by data of geometry.
We will comment upon this correspondence from categorical point of view.

We first consider a category of analytic submanifold in $\EE^n$. We refer it
$Submfd$, whose object is an analytic submanifold ${S\hookrightarrow \EE^n}$.
For two submanifolds ${S\hookrightarrow \EE^n}$ and ${S'\hookrightarrow 
\EE^n}$, there is an  analytic map $f$ between them as a morphism in the 
category {\it Submfd}. For these submanifolds,  we naturally define a 
correspondence between ${\Nabla}_{S\hookrightarrow \EE^n}^{\Dir}$ and 
${\Nabla}_{S'\hookrightarrow \EE^n}^{\Dir}$. We will introduce a category of
Dirac operator $Dirac_{\cdot \hookrightarrow \EE^n}$ whose object is the 
formal infinite series ring $\Cal C_S^{\omega\CC}[[ 
{\Nabla}_{S\hookrightarrow \EE^n}^{\Dir}]]$ and morphism $f^*$ is induced 
from $f$. Then, we can find an equivalent functor $\sigma$ between the 
category {\it Submfd} of submanifolds in $\EE^n$ and $Dirac_{\cdot 
\hookrightarrow \EE^n}$; morphism is given as $f$.

On the other hand, we can consider an analytic category {\it Anal} 
defined over
parameter space $(s^1,\cdots,s^k,t)$, whose object is
a subset of formal power series 
$\CC[[\partial_\alpha,s^\alpha,t,\partial_t]]$.
We can naturally define a forgetful functor $\nu_{Anal}$ from the
$Dirac_{\cdot \hookrightarrow \EE^n}$ to {\it Anal}.
Then it is natural to consider a $\CC$ vector space,
$$
   {\Cal A}_{\cdot\hookrightarrow 
   \EE^n}^{\Dir}[{\Nabla}_{S\hookrightarrow \EE^n}^{\Dir}] :=\{ A\in 
   \CC[[\partial_\alpha,s^\alpha]] \subset 
   \CC[[\partial_\alpha,s^\alpha,t,\partial_t]] \ |\ [\partial_t-A, 
   {\Nabla}_{S\hookrightarrow \EE^n}^{\Dir}] =0 \ \} / \Cal 
   C_S^{\omega\CC}[[ {\Nabla}_{S\hookrightarrow \EE^n}^{\Dir}]].
$$
To determine ${\Cal A}_{\cdot\hookrightarrow \EE^n}^{\Dir}$
 is commutant problem 
as that in  von-Neumann algebra. The intersection  of 
commutants ${\Cal A}_{\cdot\hookrightarrow \EE^n}^{\Dir} 
[{\Nabla}_{S\hookrightarrow \EE^n}^{\Dir}] \cap 
\Cal C_S^{\omega\CC}[[ {\Nabla}_{S\hookrightarrow \EE^n}^{\Dir}]]$ 
is related to a commutative
ring. For the cases of  curves and conformal surfaces, the intersection
is associated with algebraic curve {\it e.g.}, a hyperelliptic
curve ${\roman{Spec}}(\CC[x,y]/(y^2 - f(x))$ .
As it might, guess, be a pun of someone, 
${\roman{Spect}}({\Nabla}_{S\hookrightarrow \EE^n}^{\Dir})$ 
corresponds to ${\roman{Spec}}(\CC[x,y]/(y^2 - f(x))$.
Former is of a noncommutative ring defined over differential
geometrical manifold and is expressed by a differential operator
associated with K${}^0$-group in K-theory.
Latter is of commutative ring defined over algebraic variety
and is expressed by a commutative variable 
associated with K${}^1$-group in the K-theory [Con]. 
Using the orbits, we can introduce linear topology and classify 
${\Cal A}_{\cdot\hookrightarrow \EE^n}^{\Dir}
[{\Nabla}_{S\hookrightarrow \EE^n}^{\Dir}]$ itself.
Investigations of these sets have been already
done in the framework of soliton theory.
Geometrical interpretation of soliton theory can be regarded as 
a composite functor from {\it Anal} to {\it Submfd}. 

The morphism $f$ in {\it Submfd} induces a morphism $f^{\#}$  in {\it
Anal}. Inverse of this functor means that if we could classify the 
operators set ${\Cal A}_{S\hookrightarrow 
\EE^n}^{\Dir}[{\Nabla}_{S\hookrightarrow \EE^n}^{\Dir}]$'s, we can 
classify submanifolds in {\it Submfd}. For example, let us consider a
morphism $f^{\#}$ such that ${\Cal A}_{\cdot\hookrightarrow 
\EE^n}^{\Dir}[{\Nabla}_{S\hookrightarrow \EE^n}^{\Dir}]$ and 
${\Cal A'}_{\cdot\hookrightarrow \EE^n}^{\Dir}[
{\Nabla}_{S\hookrightarrow 
\EE^n}^{\Dir}]$ are connected by an orbit of $t$ and we denote it 
$f^{\#}_t$. Its inverse functor of $f^{\#}_t$ implies that there is a 
deformation $f_t$ of a submanifolds $S$ and $S'$. Since the Dirac 
operator  $ {\Nabla}_{S\hookrightarrow \EE^n}^{\Dir}$ connects between
them, we emphasize that the Dirac operator itself may behave 
as an interpreter like a functor.

In these cases, by using ${\Cal A}_{S\hookrightarrow \EE^n}^{\Dir}$ 
we can classify  geometrical objects or 
its moduli and can investigate an (algebraic) relation
between classified objects [Mat12, 13, 15, 17].
Moreover physically speaking, 
such classification and investigation of their relations
 can be interpreted as a quantization of submanifold geometry as
we showed in [Mat12, 13, 15, 17].

As mentioned in the introduction, one of our purposes of this study is
arthmeticalization of geometry. In the case of soliton theory,
there is an infinite dimensional Lie group acting the 
${\Cal A}_{\cdot\hookrightarrow \EE^n}^{\Dir}$ and it 
governs the structure of objects of {\it Anal}.
As Klein and Lie wished, above program might enable us
to investigate {\it Submfd}. In other words, this is 
a Galois theory of submanifold. There the Dirac operator
plays the same role as $\sqrt{-1}$ in algebraic number theory.
It implies that even though invariances in the differential geometry and
harmonic map theory are usually given by quadratic form,
I believe that we must factorize or construct a "Dirac operator"
and adjoin such a object;
I consider that such a way  is a way to unify mathematics
from discovery of $\sqrt{-1}$.  
 
From physical point of view, as I showed in [Mat12],
quantized submanifold  contains various submanifolds and is 
governed by such a large group. If the quantized submanifold is 
given first,
there naturally appears the Dirac operator and Dirac field
(as an object that defines a dimensionality of 
functional space or geometry). This might give us a deep 
philosophical question what is an origin of material or fermion.

%\newpage
\fvskip

\enddocument